\documentstyle[12pt]{article}
   \font\tenmsb=msbm10 scaled\magstep 1
   \font\sevenmsb=msbm7 scaled \magstep 1
   \font\fivemsb=msbm5 scaled \magstep 1
\newfam\msbfam
      \textfont\msbfam=\tenmsb
      \scriptfont\msbfam=\sevenmsb
      \scriptscriptfont\msbfam=\fivemsb
\def\Bbb#1{{\fam\msbfam #1}}

 \font\twelvemsb=msbm10 scaled\magstep 4
   \font\tenmsb=msbm7 scaled \magstep 4
   \font\sevenmsb=msbm5 scaled \magstep 4
\newfam\msbfamy
      \textfont\msbfamy=\twelvemsb
      \scriptfont\msbfamy=\tenmsb
      \scriptscriptfont\msbfamy=\sevenmsb

\font\tenblah=msbm10 scaled\magstep 3
   \font\sevenblah=msbm7 scaled \magstep 1
   \font\fiveblah=msbm5 scaled \magstep 1
\newfam\blahfam
      \textfont\blahfam=\tenblah
      \scriptfont\blahfam=\sevenblah
      \scriptscriptfont\blahfam=\fiveblah

\font\tengothic=eufm10 scaled\magstep 1
\font\sevengothic=eufm7 scaled\magstep 1
\newfam\gothicfam
      \textfont\gothicfam=\tengothic
      \scriptfont\gothicfam=\sevengothic
\def\goth#1{{\fam\gothicfam #1}}

\newcommand{\qed}{\hskip 1cm $\rlap{$\sqcap$}\sqcup$}

\newcommand{\proj}[1]
{ \mathchoice
           { {\Bbb P}^{#1} }
           { {\Bbb P}^{#1} }
           { {\Bbb P}^{#1} }
           { {\Bbb P}^{#1} }
         }

\newcommand{\zed}{ {\Bbb Z} }

\newcommand{\verylong}{\hbox{$\hbox to .45in{\rightarrowfill}$}  }
\newcommand{\kindalong}{\hbox{$\hbox to .35in{\rightarrowfill}$}  }

\newtheorem{thm}{Theorem}[section]
\newtheorem{prop}[thm]{Proposition}
\newtheorem{lemma}[thm]{Lemma}
\newtheorem{cor}[thm]{Corollary}

\newtheorem{defn}[thm]{Definition}
\newtheorem{rmk}[thm]{Remark}
\newtheorem{exemp}[thm]{Example}

\newtheorem{conj}[thm]{Conjecture}
\newtheorem{fff}[thm]{Fact}

\newenvironment{example}{\begin{exemp}\begin{rm}}{\end{rm} \qed
\end{exemp}}
\newenvironment{remark}{\begin{rmk}\begin{rm}}{\end{rm} \qed
\end{rmk}}
\newenvironment{definition}{\begin{defn}\begin{rm}}{\end{rm} \qed
\end{defn}}
\newenvironment{fact}{\begin{fff}\begin{rm}}{\end{rm} \qed
\end{fff}}

\begin{document}
\title{Determinantal Schemes and Buchsbaum-Rim Sheaves}

\author{Martin Kreuzer \\
        Fakult\"at f\"ur Mathematik \\
        Universit\"at Regensburg \\
        D-93040 Regensburg \\
        Germany
\and
Juan C. Migliore \\
        Department of Mathematics \\
        University of Notre Dame \\
        Notre Dame, IN 46556 \\
        USA
\and
Uwe Nagel \\
FB Mathematik-Informatik-17 \\
Universit\"at-GH Paderborn \\
D-33 095 Paderborn \\
Germany
   \and
        Chris Peterson \\
        Department of Mathematics \\
        University of Notre Dame \\
        Notre Dame, IN 46556 \\
        USA }
\date{ }
\maketitle
\footnoterule{{\em 1991 Mathematics Subject Classification:} Primary 14M12,
14F05; Secondary 13D02, 14C20, 13C40}

\baselineskip 15pt



\section{Introduction}

A natural and efficient method for producing numerous examples of interesting
schemes is to consider the vanishing locus of the minors of a homogeneous
polynomial matrix. If the matrix satisfies certain genericity conditions then
the resulting schemes have a number of well described properties. These objects
have been studied in both a classical context and a modern context and go
by the
name of determinantal schemes. Some of the classical schemes that can be
constructed in this manner are the Segre varieties, the rational normal
scrolls,
and the Veronese varieties. In fact, it can be shown (cf.\ \cite{harris}) that
any projective variety is isomorphic to a determinantal variety arising from a
matrix with linear entries! Due to their important role in algebraic
geometry and
commutative algebra, determinantal schemes and their associated rings have
both merited and received considerable attention in the literature.
Groundbreaking work has been carried out by a number of different authors; we
direct the reader to the two excellent sources \cite{bruns-vetter} and \cite
{eisenbud} for background, history, and a list of important
papers.

 A homogeneous polynomial matrix can be viewed as defining a map between
free modules defined over the underlying polynomial ring. Associated to
such a map are a number of complexes. The most important of these
are the Eagon-Northcott and Buchsbaum-Rim complexes. Under appropriate
genericity conditions, these complexes are exact and it is in this
special situation where we will focus our attention. Buchsbaum-Rim sheaves
are a
family of sheaves associated to the sheafified Buchsbaum-Rim  complex.  In
particular, a first Buchsbaum-Rim sheaf is the kernel of a generically
surjective map between two direct sums of line bundles, whose cokernel is
supported in the correct codimension.  This family of sheaves is described and
studied in the two papers \cite{mig-pet}, \cite{MNP}.

A certain aspect of these sheaves was found to bear an interesting relationship
to earlier work of the first author. In \cite{kreuzer}, Kreuzer obtained the
following  characterization of 0-dimensional complete intersections in
$\proj{3}$:

\bigskip

\noindent {\bf Theorem} (\cite{kreuzer} Theorem 1.3) \ \ {\em A 0-dimensional
subscheme $Y \subset \proj{3}$ is a complete intersection if and only if $Y$ is
arithmetically Gorenstein and there exists an arithmetically Cohen-Macaulay,
l.c.i.\ curve $C$ such that $Y$ is the associated subscheme of an effective
Cartier divisor on $C$ and ${\cal O}_C (Y) \cong \omega_C (-a_Y )$ is globally
generated.}

\bigskip

Complete intersections form a very important subset of the more general class
of standard determinantal schemes (i.e the determinantal subschemes of
$\proj{n}$ arising from the maximal minors of a homogeneous matrix of the
``right size"). One immediately observes that to every standard determinantal
scheme is associated a number of Buchsbaum-Rim sheaves and to every
Buchsbaum-Rim sheaf is associated a standard determinantal ideal. We say a
standard determinantal scheme is ``good" if one can delete a generalized row
from its corresponding matrix and have the maximal minors of the resulting
submatrix define a scheme of the expected codimension. In particular, complete
intersections are good, as are most standard determinantal schemes.

The paper is
organized as follows. In Section 2 we provide the necessary  background
information. The next section is the heart of the paper. Here we give several
characterizations of standard and good determinantal subschemes. Some  of these
results are summarized in the following:

\bigskip

\noindent {\bf Theorem}{\em \ \
Let $X$ be a subscheme of $\proj{n}$ with $codim \ X
\geq 2$. The following are equivalent.
\newcounter{tempA}
\begin{list} {(\alph{tempA})}{\usecounter{tempA}}
    \setlength{\rightmargin}{\leftmargin}
\item $X$ is a good determinantal scheme of codimension $r+1$.

\item $X$ is the zero-locus of a regular section of the dual of a first
Buchsbaum-Rim sheaf of rank $r+1$.

\item $X$ is standard determinantal and locally a complete intersection outside
a subscheme $Y \subset X$ of codimension $r+2$ in $\proj{n}$.

\end{list}
}
\noindent Several of our results in Section~\ref{char-good} involve the
cokernel of the map of free modules mentioned above.  We do not quote
these results here since we need some notation from Section~\ref{prelim-sect}.
These results are  important in Section~\ref{corollaries}, though, where we
give
our main generalizations of Kreuzer's theorem.  We mention two of these.

\bigskip

\noindent {\bf Corollary}{\em \ \
Let $X \subset \proj{n}$ be a subscheme of codimension $r+1 \geq 3$. Then
$X$ is a complete intersection if and only if $X$ is arithmetically
Gorenstein and there is a good determinantal subscheme $S \subset \proj{n}$
of codimension $r$ and a canonically defined sheaf ${\cal M}_S$ on $S$ (in
codimension two, ${\cal M}_S \cong \omega_S$ up to twist)
such
that $X \subset S$ is the zero-locus of a regular section  $t \in H^0_* (S,
{\cal M}_S )$. Furthermore, $S$ and
${\cal M}_S$ can be chosen so that ${\cal M}_S$ is globally generated.
}

\bigskip

\noindent {\bf Corollary}{\em \ \
Suppose $X \subset \proj{3}$ is zero-dimensional.  Then the following are
equivalent:
\newcounter{tempB}
\begin{list} {(\alph{tempB})}{\usecounter{tempB}}
    \setlength{\rightmargin}{\leftmargin}
\item $X$ is good determinantal;

\item $X$ is standard determinantal and a local complete intersection;

\item There is an arithmetically Cohen-Macaulay curve $S$, which is a local
  complete
intersection, such that $X$ is a subcanonical Cartier divisor on $S$.

\end{list}
Furthermore, $X$ is defined by
a $t \times (t+r)$ matrix if and only if the Cohen-Macaulay type of $X$ is
${{r+t-1} \choose r}$ and that of
$S$ is ${{r+t-1}
\choose {r-1}}$. }

\bigskip

\noindent The last sentence of this corollary gives the connection to
Kreuzer's
theorem: recall that the only standard determinantal subschemes with
Cohen-Macaulay type 1 (i.e.\ arithmetically Gorenstein) are complete
intersections.
 In a similar way we characterize  good determinantal subschemes
of $\proj{n}$ of any codimension, with special, stronger, results in the
case of
zeroschemes and the case of codimension two subschemes.   We close with a
number
of examples.


\section{Preliminaries}\label{prelim-sect}

Let $R=k[x_0,x_1,\dots,x_n]$ be a  polynomial ring with the standard grading,
where
$k$ is an infinite field and $n\geq 2$.
For any sheaf $\cal F$ on $\proj{n}$, we define
$H^i_*(\proj{n},{\cal F})=\bigoplus_{t\in \zed}H^i(\proj{n},{\cal F}(t))$.
For any scheme $V \subset \proj{n}$, $I_V$ denotes the saturated homogeneous
ideal of $V$ and ${\cal I}_V$ denotes the ideal sheaf of
$V$ (hence $I_V= H^0_* (\proj{n},{\cal I}_V)$).

\begin{definition}
If $A$ is a homogeneous matrix, we denote by $I(A)$ the ideal of maximal minors
of $A$.  A codimension $r+1$ scheme, $X$, in $\proj{n}=Proj(R)$ will be
called a
{\em standard determinantal scheme} if $I_X = I(A)$ for some homogeneous $t
\times (t+r)$ matrix, $A$. $X$ will be called a {\em good determinantal scheme}
if additionally, $A$ contains a $(t-1) \times (t+r)$ submatrix (allowing a
change
of basis if necessary-- see Example~\ref{good-ex}) whose ideal of maximal
minors
defines a scheme of codimension $r+2$.  In a similar way we define standard and
good determinantal ideals.
\end{definition}

\begin{example}
The ideal defined by the maximal minors of the matrix
\[
\left [
\begin{array}{cccccccccc}
x_1 & x_2 & x_3 & 0 \\
0 & x_1 & x_2 & x_3
\end{array}
\right ]
\]
is an example of a standard determinantal ideal which is not good.  Note that
this ideal is the square of the ideal of a point in $\proj{3}$, and is not a
local complete intersection (see Proposition~\ref{good-det}).
\end{example}

Note that standard determinantal schemes form an important subclass of
the more general notion of determinantal schemes, where smaller minors are
allowed (among other generalizations).  See for
instance \cite{bruns-vetter}, \cite{eisenbud}, \cite{harris}.

\begin{remark}\label{one-minor}
In the next section we will make a deeper study of good determinantal
schemes.  For now, though, we observe the following.  Let $X$ be a standard
determinantal scheme coming from a $t \times (t+r)$ matrix $A$.  Then $X$ is
good if and only if there is a $(t-1) \times (t-1)$ minor of $A$ which does not
vanish on any component of $X$ (possibly after making a change of basis).  In
particular, we formally include the possibility that $t=1$, and we include the
complete intersections among the good determinantal schemes.
\end{remark}

\begin{fact} \label{EN-and-BR}
Let ${\cal F}$ and ${\cal G}$ be locally free sheaves of ranks $f$ and $g$
respectively on a smooth variety $Y$. Let $\phi :{\cal F}\rightarrow {\cal G}$
be a generically surjective homomorphism. We can associate to $\phi$ an
Eagon-Northcott complex

\begin{equation}\label{ENSeq}
\begin{array}{cccc}
0 \rightarrow
\wedge^f{\cal F} \otimes (S^{f-g}
{\cal G})^{\vee} \otimes \wedge^g {\cal G}^{\vee}
\rightarrow
\wedge^{f-1} {\cal F} \otimes (S^{f-g-1} {\cal G})^{\vee} \otimes \wedge^g
{\cal G}^{\vee}
\rightarrow\dots \\
\hskip1.5cm \rightarrow
\wedge^{g+1} {\cal F} \otimes {\cal G}^{\vee} \otimes \wedge^g
{\cal G}^{\vee}
\rightarrow
\wedge^g {\cal F} \otimes \wedge^g {\cal G}^{\vee}
\buildrel \wedge^g \phi \over \rightarrow {\cal O}_Y \rightarrow 0\\
\end{array}
\end{equation}
and a Buchsbaum-Rim complex
\begin{equation}\label{BRSeq}
\begin{array}{cccc}
0 \rightarrow \wedge^f {\cal F} \otimes S^{f-g-1} {\cal G}^\vee \otimes
\wedge^g {\cal G}^\vee \rightarrow \wedge^{f-1} {\cal F} \otimes S^{f-g-2}
{\cal G}^\vee \otimes \wedge^g {\cal G}^\vee \rightarrow \dots \\
\hskip 1cm \rightarrow \wedge^{g+2} {\cal F} \otimes {\cal G}^\vee \otimes
\wedge^g {\cal G}^\vee \rightarrow \wedge^{g+1} {\cal F} \otimes \wedge^g
{\cal G}^\vee \rightarrow {\cal F} \buildrel \phi \over \rightarrow {\cal G}
\rightarrow 0\\
\end{array}
\end{equation}
(see \cite{GLP}, \cite{eisenbud}, \cite{BR}, \cite{EN}, \cite{buchs64}).    If
the support of the cokernel of $\phi$ has the expected codimension $f-g+1$
then these complexes are acyclic.
\end{fact}
The consequences of this fact will play a
crucial role throughout the paper and they lead us to the following definition.

\begin{definition}\label{def-BR}
Let ${\cal F}$ and ${\cal G}$ be two locally
free sheaves which split as the sum of line bundles and let $\phi :{\cal
F}\rightarrow {\cal G}$ be a generically surjective homomorphism whose
cokernel is supported on a scheme with the ``expected'' codimension $f-g+1$. As
mentioned in the fact above, the Buchsbaum-Rim complex will be exact and
provides a free resolution of the cokernel of the map $\phi$. The kernel of the
map $\phi$ will be called a {\it first Buchsbaum-Rim sheaf}.
We use the symbol ${\cal B}_{\phi}$ to represent such a sheaf.
\end{definition}

More generally, the $i^{th}$ Buchsbaum-Rim sheaf associated to $\phi$ is the
$(i+1)^{st}$ syzygy sheaf in the Buchsbaum-Rim complex.  However, in
this paper we will use only the first Buchsbaum-Rim sheaves.

\begin{remark} \label{free-is-BR}
In Fact~\ref{EN-and-BR} and Definition~\ref{def-BR}, we will allow the rank of
$\cal G$ to be zero, and use the convention that even in this case, $\wedge^0
{\cal G}^\vee = {\cal O}_Y$.  Moreover, the Buchsbaum-Rim complex
becomes $0 \rightarrow {\cal F} \rightarrow {\cal F} \buildrel \phi \over
\rightarrow 0$, and it follows that the sheafification of any free module is a
first Buchsbaum-Rim sheaf.

In Fact~\ref{EN-and-BR} and Definition~\ref{def-BR}, we can also start with
free modules $F$ and $G$, and we get Eagon-Northcott and Buchsbaum-Rim
complexes
of free modules.  The corresponding kernel of the map $\phi$ will then be
called a {\em first Buchsbaum-Rim module}.  Note that in this context $\phi$
can be represented by a homogeneous matrix $\Phi$, and the image of $\wedge^g
\phi$ is precisely $I(\Phi )$.

Note also that since first Buchsbaum-Rim sheaves (resp.\ modules) are
second syzygy sheaves (resp.\ modules), they are reflexive.
\end{remark}

\begin{fact} (\cite{eisenbud} exer.\ 20.6 or \cite{BE77})\label{annihilator}
Let $\Phi$ be a matrix whose ideal $I(\Phi )$ of maximal minors vanishes in the
expected codimension, and so $coker \ \Phi$ has a corresponding Buchsbaum-Rim
resolution.  Then the annihilator of $coker \ \Phi$ is precisely $I(\Phi )$.
\end{fact}

In this paper, we will often be interested in going in the opposite direction,
starting with a standard determinantal ideal $J$ and considering the possible
associated matrices and cokernels.  With this in mind, we make the following
definition.

\begin{defn}\label{Mv} \begin{rm}
Let $X$ be a standard determinantal scheme of codimension $r+1$ with
corresponding ideal $I_X$.  Then we set
\[
{\cal M}_X := \left \{
\begin{array}{c|c}
M &
\begin{array}{l}
\hbox{\rm $M$ is a f.g.\ graded $R$-module with $Ann_R M = I_X$ and a minimal}
\\
\hbox{\rm presentation of the form $\displaystyle R^{r+\mu} \rightarrow R^\mu
\rightarrow M \rightarrow 0$}
\end{array}

\end{array}
\right \}
\]
\end{rm}
\end{defn}

\bigskip

\noindent ${\cal M}_X$ is the set of possible cokernels of homogeneous matrices
whose ideals of maximal minors are precisely $I_X$.  In some situations, ${\cal
M}_X$ consists of just one element (up to isomorphism and twisting).  For
example, it can be shown that this happens if $r=1$ (i.e.\ codimension 2, using
Hilbert-Burch theory-- see Corollary~\ref{codim2}).  ${\cal M}_X$ also consists
of just one element if
$X$ is a complete intersection.  We do not know the precise conditions which
guarantee that all the elements of ${\cal M}_X$ are
isomorphic up to twisting.  In any case, we can at least show that the
elements of ${\cal M}_X$ look very much alike:

\begin{lemma}
The elements of ${\cal M}_X$ all have the same graded Betti numbers, up to
twisting, and in particular come from matrices of the same size.
\end{lemma}

\noindent {\em Proof:}
Let $M_1, M_2 \in {\cal M}_X$ and assume that $M_i$ has $t_i$ minimal
generators, $i = 1,2$.  We may also assume that $M_i$ is the cokernel of
a $t_i \times (t_i +r)$ matrix
$\Phi_i$.  By \cite{eisenbud} p.\ 494, $Rad(I(\Phi )) = Rad (Ann_R M_i) = Rad
(I_X )$.  Hence $I(\Phi)$ is a homogeneous matrix  defining a subscheme of
$\proj{n}$ of codimension $r+1$, the expected codimension, and we may apply the
Eagon-Northcott complex to get a minimal free resolution for $I(\Phi ) = I_X$.
Hence $I_X$ has ${{r+t_1} \choose r} = {{r+t_2} \choose r}$ minimal generators,
and $t_1 = t_2$.

Now let $M \in {\cal M}_X$ and assume that it has $t$ minimal  generators.
There is a minimal free resolution
\[
\cdots \rightarrow F \buildrel \Phi \over \rightarrow G \rightarrow M
\rightarrow 0
\]
where $rk \ F = t+r$ and $rk \ G = t$.  As above, $I(\Phi)$ defines a
subscheme of codimension $r+1$, and so the Buchsbaum-Rim complex resolves
$M$ and we are done. \qed

\begin{prop}\label{entries-of-syz}
Let ${\cal F}$ and ${\cal G}$ be locally free sheaves of ranks $f$ and
$g$ respectively on $\proj{n}$. Let $\phi :{\cal
F}\rightarrow {\cal G}$ be a generically surjective homomorphism. Assume the
cokernel of $\phi$ is supported on a scheme of codimension $f-g+1$. Let
$I_{\phi}$ denote the homogeneous ideal of the scheme determined by the
cokernel
of
$\wedge^g\phi$. Let $I_s$ denote the homogeneous ideal of the
zero-locus of a section, $s\in H^0(\proj{n},{\cal B}_{\phi})$ (where
${\cal B}_{\phi}$ denotes the local first Buchsbaum-Rim sheaf of $\phi$).
Let $I_t$ denote the homogeneous ideal of the
zero-locus of a section, $t\in H^0(\proj{n},{\cal B}_{\phi}^*)$ (where
${\cal B}_{\phi}^*$ denotes the dual of ${\cal B}_{\phi}$).
Then for any such section, $I_s\subset I_{\phi}$ and $I_t\subset I_{\phi}$
\end{prop}

\noindent {\em Proof:}
Locally, we can represent the map $\phi$ by an $f\times g$ matrix, $A$. In the
same local coordinates, the map from $\wedge^{g+1}{\cal F}\otimes \wedge^g{\cal
G}^{\vee}$ to ${\cal F}$ (in the Buchsbaum-Rim complex associated to $\phi$)
can be expressed by a matrix, $M$. The entries of $M$ can be written in terms
of $A$ as follows. Let $I_A$ denote the ideal of maximal minors of the matrix
$A$. $I_A$ locally describes the scheme defined by $I_{\phi}$. Each column in
the matrix, $M$, arises from choosing $t+1$ columns of the matrix $A$ and
considering all $t\times t$ minors of this submatrix of $A$. Thus, each entry
in the matrix $M$ is an element of $I_A$. Locally, sections of the first
Buchsbaum-Rim sheaf of $\phi$ are determined by an element of the column
space of $M$ (considered as a module). An immediate consequence of this fact is
that the vanishing locus of any section of the first Buchsbaum-Rim sheaf of
$\phi$ or the dual of the first Buchsbaum-Rim sheaf of $\phi$ will contain the
scheme defined by
$I_{\phi}$. \qed

\begin{remark}
For clarity, and because of its importance, we restrict ourselves to
determinantal  subschemes of projective space in the body of this paper.
However, the reader will observe that many of our arguments hold true for
subschemes of a  smooth projective variety and some even for determinantal
ideals
of an arbitrary  commutative ring.
\end{remark}


\section{Characterizations of Good Determinantal Schemes}\label{char-good}

In \cite{mig-pet} and \cite{MNP}, regular sections of first Buchsbaum-Rim
sheaves were considered, and it was shown that they possess many interesting
properties.  For example, a regular section of a first Buchsbaum-Rim sheaf of
odd rank has a zero-locus whose top dimensional part is arithmetically
Gorenstein.

In this paper we are primarily concerned with regular sections of the {\em
dual} of a first Buchsbaum-Rim sheaf.  Our first result gives a property
which is
analogous to the ones mentioned above for the first Buchsbaum-Rim sheaves.

\begin{thm}\label{good-iff-sect}
Let $X$ be a subscheme of $\proj{n}$ with $codim \ X
\geq 2$. The following are equivalent.
\newcounter{temp}
\begin{list} {(\alph{temp})}{\usecounter{temp}}
    \setlength{\rightmargin}{\leftmargin}
\item $X$ is a good determinantal scheme of codimension $r+1$.

\item $X$ is the zero-locus of a regular section of the dual of a first
Buchsbaum-Rim sheaf of rank $r+1$.

\end{list}
\end{thm}

\noindent {\em Proof:}
We first prove (a) $\Rightarrow$ (b).  By assumption there is a homomorphism
$\Phi$ such that  $I_X = I(\Phi )$, and we have an exact sequence
\begin{equation}
0 \rightarrow B \rightarrow F \buildrel \Phi \over \rightarrow G \rightarrow
coker \ \Phi \rightarrow 0
\end{equation}
where $rk \ G = t$, $rk \ F = t+r$ and $B$ is a first Buchsbaum-Rim module.

If $t=1$ then $I(\Phi )$ is a complete intersection of height $r+1$, which can
be viewed as a section of (the dual of) a free module of rank $r+1$.  By
Remark~\ref{free-is-BR}, a free module is a first Buchsbaum-Rim module.
Hence we
can assume from now on that $t \geq 2$.

Since $X$ is a good determinantal scheme, there is a projection $\pi : G
\rightarrow G'$, where $G'$ has rank $t-1$, $G'$ is obtained from $G$ by
removing
one free summand $R(a)$, and such that $ht ( I(\pi \circ \Phi )) = r+2$.
We get
a commutative diagram
\begin{equation}\label{usual-diag}
\begin{array}{ccccccccccc}
&&&&&& 0 \\
&&&&&& \downarrow \\
&&&& 0 & \rightarrow & R(a) & \rightarrow & R(a) & \rightarrow & 0 \\
&&&& \downarrow && \downarrow \\
0 & \rightarrow & B & \rightarrow & F & \buildrel \Phi \over \rightarrow & G &
\rightarrow & coker \ \Phi & \rightarrow & 0 \\
&&&& || && \phantom \pi \downarrow \pi \\
0 & \rightarrow & B' & \rightarrow & F & \buildrel {\Phi '} \over \rightarrow &
G' & \rightarrow & coker \ \Phi ' & \rightarrow & 0 \\
&&&& \downarrow && \downarrow \\
&&&& 0 && 0 \\

\end{array}
\end{equation}

Let $\alpha$ be the induced injection from $B$ to $B'$.  Twist everything in
(\ref{usual-diag}) by $-a$ and relabel, so that the Snake Lemma gives that
$I =
coker \ \alpha$ is an ideal and we have an exact sequence
\begin{equation}\label{short-exact}
0 \rightarrow R/I \rightarrow coker \ \Phi \rightarrow coker \ \Phi '
\rightarrow 0
\end{equation}
It follows that $I_X = I(\Phi ) = Ann ( coker \ \Phi) \subset I$ (see
Fact~\ref{annihilator}), where $I_X$ is the saturated ideal of $X$.

On the other hand, it follows from the same exact sequence that
\[
Ann ( coker \ \Phi' ) \cdot I \subset Ann(coker \ \Phi ) = I(\Phi ) = I_X.
\]
But since $X$ is good determinantal, it follows that $I(\Phi' ) = Ann ( coker \
\Phi' )$ and
$ht ( I(\Phi' )) > ht (I(\Phi ))$.  Hence $I \subset I(\Phi )$ and so we
conclude $I = I(\Phi ) = I_X$.  But then we have a short exact sequence
\[
0 \rightarrow B \rightarrow B' \rightarrow I_X \rightarrow 0
\]
and so by sheafifying, it follows that $X$ is the zero-locus of a regular
section of the dual of the first Buchsbaum-Rim sheaf ${\cal B}'$ as claimed.
(Note that $B'$ is reflexive-- see Remark~\ref{free-is-BR}.)

We now prove (b) $\Rightarrow$ (a).  Assume that $X$ is the zero-locus of a
regular section of a sheaf $({\cal B}')^*$, where ${\cal B'}$ is the
sheafification of a first Buchsbaum-Rim module $B'$ of rank $r+1$. We are thus
given exact sequences (after possibly replacing $B'$ by a suitable twist)
\begin{equation}\label{BRseq}
0 \rightarrow B' \rightarrow F \buildrel {\Phi'} \over \longrightarrow G
\rightarrow coker \ \Phi' \rightarrow 0
\end{equation}
and
\begin{equation}\label{sect}
0 \rightarrow R \rightarrow (B')^* \rightarrow Q \rightarrow 0
\end{equation}
such that $rk \ F = t+r$, $rk \ G = t-1$, $Ann ( coker \ \Phi' ) =
I(\Phi ')$ (which has height $r+2$) and
\[
0 \rightarrow Q^* \rightarrow B' \rightarrow I \rightarrow 0
\]
is exact (again, $B'$ is reflexive), where $I$ is an ideal whose saturation is
$I_X$.  One can check that dualizing (\ref{BRseq}) provides
\[
0 \rightarrow G^* \rightarrow F^* \rightarrow (B')^* \rightarrow 0.
\]
The mapping cone procedure applied to (\ref{sect}) then gives
\[
0 \rightarrow R \oplus G^* \rightarrow F^* \rightarrow Q
\rightarrow 0.
\]
Dualizing this, we obtain the following commutative diagram:
\[
\begin{array}{cccccccccccc}
&&&&&& 0 \\
&&&&&& \downarrow \\
&&0 && 0 && R  \\
&& \downarrow && \downarrow && \downarrow \\
0 & \rightarrow & Q^* & \rightarrow & F & \buildrel \Phi \over \rightarrow
  & R \oplus G & \rightarrow & coker \ \Phi & \rightarrow & 0 \\
&& \downarrow && || && \downarrow && \downarrow \\
0 & \rightarrow & B' & \rightarrow & F & \buildrel {\Phi '} \over \rightarrow
& G & \rightarrow & coker \ \Phi ' & \rightarrow & 0 \\
&& \downarrow && \downarrow && \downarrow \\
&& I && 0 && 0 \\
&& \downarrow \\
&& 0
\end{array}
\]
The Snake Lemma then gives
\[
0 \rightarrow R/I \rightarrow coker \ \Phi \rightarrow coker \ \Phi'
\rightarrow 0.
\]
It follows that
\[
I \cdot Ann (coker \ \Phi' ) = I \cdot I(\Phi' ) \subset Ann(coker \ \Phi ).
\]
Thus $ht (Ann ( coker \ \Phi )) \geq r+1$.  Note that the maximal possible
height of $Ann ( coker \ \Phi )$ is $r+1$, hence we get $ht (Ann ( coker \ \Phi
)) = r+1$ and  $Q^*$ is a first Buchsbaum-Rim module.  From the Buchsbaum-Rim
complex one can then check that
$H^1_* (\proj{n} , {\cal Q}^* ) = 0$, and hence $I = I_X$ is saturated.
Then as
in the first part we get
$I_X = I(\Phi )$, as desired.
  \qed

\medskip

We now give a result which characterizes the good determinantal schemes among
the standard determinantal schemes.  We use the set ${\cal M}_X$ introduced in
Definition~\ref{Mv}.

\begin{prop}\label{good-det}
Suppose that $X$ is a standard determinantal scheme of codimension $r+1$.  Then
the following are equivalent.
\newcounter{temp2}
\begin{list} {(\alph{temp2})}{\usecounter{temp2}}
    \setlength{\rightmargin}{\leftmargin}
\item $X$ is good determinantal;

\item There is an $M_X \in {\cal M}_X$ and an embedding $R/I_X \hookrightarrow
M_X$ whose image is a minimal generator of $M_X$ as an $R$-module, and whose
cokernel is supported on a subscheme of codimension $\geq r+2$.

\item There is an element $M_X \in {\cal M}_X$ which is an ideal in $R/I_X$ of
positive height.

\end{list}

Furthermore, if any of the above conditions hold then $X$ is a local complete
intersection outside a subscheme $Y \subset \proj{n}$ of codimension $r+2$.

\end{prop}

\begin{remark} The first two parts of the above proposition do not even require
that the field be infinite.
\end{remark}

\noindent {\em Proof of \ref{good-det}}
We begin with (a) $\Rightarrow$ (b).  Assume that $X$ is a good determinantal
scheme arising from a homogeneous matrix $\Phi$.  As in the proof of
Theorem~\ref{good-iff-sect} (see the diagram (\ref{usual-diag})), we have
(after
possibly twisting) a commutative diagram
\begin{equation}
\begin{array}{ccccccccccc}
&&&&&&  && 0 \\
&&&&&&  && \downarrow  \\
&& 0 &&  &  & R & \rightarrow & R/I_X &  \\
&& \downarrow &&  && \downarrow && \downarrow \\
0 & \rightarrow & B & \rightarrow & F & \buildrel \Phi \over \rightarrow & G &
\rightarrow & M_X & \rightarrow & 0 \\
&& \downarrow && || && \phantom \pi \downarrow \pi  && \downarrow \\
0 & \rightarrow & B' & \rightarrow & F & \buildrel {\Phi '} \over \rightarrow &
G' & \rightarrow & M_Y & \rightarrow & 0 \\
&&\downarrow &&  && \downarrow && \downarrow \\
&& I_X &&&& 0 && 0 \\
&& \downarrow \\
&& 0

\end{array}
\end{equation}
where $rk \ F = t+r$, $rk \ G = t$, $rk \ G' = t-1$, $\Phi'$ is obtained by
deleting a suitable row of
$\Phi$,
$Y$ is the codimension $r+2$ scheme defined by the maximal minors of $\Phi'$,
$B$ and $B'$ are the kernels of $\Phi$ and $\Phi'$, respectively, and $M_X$ and
$M_Y$ are the respective cokernels.  Then all parts of (b) follow immediately.

This diagram also proves the last part of the Proposition, since by
Theorem~\ref{good-iff-sect}  $X$ is the zero-locus of a section of ${\cal
B}'$, the sheafification of $B'$, which is locally free of rank $r+1$ outside
$Y$.

We now prove (b) $\Rightarrow$ (a).  The assumptions in (b) imply a commutative
diagram
\[
\begin{array}{cccccccccc}
&& 0 && 0 \\
&& \downarrow && \downarrow \\
&& R & \rightarrow & R/I_X & \rightarrow & 0\\
&& \downarrow && \phantom s \downarrow s \\
F & \buildrel \Phi \over \rightarrow & G & \rightarrow & M_X & \rightarrow
& 0\\
&& \phantom \alpha \downarrow \alpha && \downarrow \\
&& G' & \buildrel \beta \over \rightarrow & coker \ s & \rightarrow & 0 \\
&& \downarrow && \downarrow \\
&& 0 && 0

\end{array}
\]
with $rk \ F = t+r$, $rk\ G = t$, $rk \ G' = t-1$.  Define $\Phi' = \alpha
\circ
\Phi$.  One can then show that
\[
F \buildrel {\Phi '} \over \rightarrow G' \buildrel \beta \over \rightarrow
coker \ s \rightarrow 0
\]
is exact.  (Either use a mapping cone argument, splitting off $R$, or else use
a somewhat tedious diagram chase.)  The assumption on the support of the
cokernel of $s$ implies $height (I(\Phi')) = r+2$, so $X$ is good, proving (a).

Now we prove (a) $\Rightarrow $ (c).  The assumption that $X$ is good
implies, in particular, that the ideal of $(t-1) \times (t-1)$ minors of $\Phi$
has height $\geq r+2$.  Hence after possibly making a change of basis, we can
apply Remark~\ref{one-minor} and \cite{eisenbud} Theorem~A2.14 (p.\ 600) to
obtain $M_X = coker \ \Phi \cong J/I_X$, where $J \subset R$ is an ideal of
height $\geq r+2$, proving (c).

Finally we prove (c) $\Rightarrow$ (b).  Since $M_X$ is an ideal of positive
height in $R/I_X$, we can find $f \in R$ with $\bar f = f \ mod \ I_X \in M_X$
such that the map $R/I_X \buildrel s \over \rightarrow M_X , \ 1 \mapsto
\bar f$
is injective.  We can even choose $f$ so that $\bar f$ is a minimal generator
of $M_X$, considered as an $R$-module.  Then $coker \ s \cong M_X / (\bar f
\cdot R/I_X )$ shows that $I_X + (f) \subset Ann_R (coker \ s)$, so $coker \ s$
is supported on a subscheme of height $\geq r+2$.  \qed
\bigskip

Next, we want to give an intrinsic characterization of good determinantal
subschemes.

\begin{thm}\label{good-iff-glci}
Suppose that $codim \ X = r+1$.  Then the following are equivalent:
\newcounter{temp5}
\begin{list} {(\alph{temp5})}{\usecounter{temp5}}
    \setlength{\rightmargin}{\leftmargin}
\item $X$ is good determinantal;

\item $X$ is standard determinantal and locally a complete intersection outside
a subscheme $Y \subset X$ of codimension $r+2$ in $\proj{n}$.

\end{list}
\end{thm}

\noindent {\em Proof:}
In view of Proposition~\ref{good-det}, we only have to prove (b) $\Rightarrow$
(a).  We again start with the exact sequence
\[
0 \rightarrow B \rightarrow F \buildrel \Phi \over \rightarrow G \rightarrow
M_X \rightarrow 0
\]
where $F$ and $G$ are free of rank $t+r$ and $t$ respectively.

Now let $P$ be a point of $X$ outside $Y$, with ideal $\wp \subset R$.  By
assumption, $X$ is a complete intersection at $P$.  We first claim that
$(M_X)_\wp \cong (R/I_X )_\wp$.  To see this, we first note that
localizing
$\Phi$ at $\wp$, we  can split off, say, $s$  direct summands until  the
resulting  map is minimal. Then the ideal of maximal minors of this matrix  has
precisely ${{r+t-s}
\choose {t-s}}$ minimal generators (Eagon-Northcott complex). On the other hand
it is a complete intersection, hence $t-s = 1$ and the cokernel $(M_X )_\wp$ of
$\Phi_\wp$ is as claimed.

Using the above isomorphism, we note that $(M_X )_\wp$ has exactly one  minimal
generator as an $R_\wp$-module.  Then by \cite{bruns-vetter},
Proposition~16.3, it follows that the ideal of submaximal minors of $\Phi$ is
not contained in $\wp$.  Since $P$ was chosen to be any point outside of
$Y$ and $codim \ Y = r+2$, it follows that no component of $X$ lies in the
ideal
of submaximal minors.  That is, the ideal of submaximal minors has height
greater than that of $I_X$.  Hence by \cite{eisenbud} p.\ 600, Theorem A2.14,
we can conclude that $M_X$ is an ideal in $R/I_X$ of positive height.
Therefore $X$ is good determinantal, by Proposition~\ref{good-det}, (c). \qed

\begin{rmk}\label{gci}\begin{rm}
Recall that a subscheme of $\proj{n}$ is said to be a {\em generic complete
intersection}
if it is locally a complete intersection at all its components. In
particular, every integral subscheme is a generic complete
intersection. This notion occurs naturally in the Serre correspondence
which relates reflexive sheaves and generic complete intersections of
codimension two (cf., for example, \cite{H2}).

Since the locus of points at which a subscheme fails to be locally a
complete intersection is closed, for a subscheme $X$ of codimension $r+1$ the
conditions being a generic complete intersection and being  locally a
complete intersection outside a subscheme $Y \subset X$ of codimension $r+2$
in $\proj{n}$ are equivalent. Thus we can reformulate the last result as
follows:
\begin{quote}
A subscheme is good determinantal if and only if it is standard
determinantal and a generic complete intersection.  \qed
\end{quote}  \end{rm}
\end{rmk}

\begin{lemma}\label{alg-lemma}
Let $A$ be a ring and let ${\goth a} \subset A$ be an ideal containing an
$A$-regular element $f$.  Let ${\goth b} := fA :_A {\goth a} = Ann_A ({\goth a}
/fA)$.  Then
$Hom_A ({\goth a}, A) \cong {\goth b}$.
\end{lemma}

\noindent {\em Proof:}
If $grade \ {\goth a} \geq 2$ then it is well-known that  $Hom_A ({\goth a}, A)
\cong A$ (up to shift in the graded case).  The interesting
case is $grade \ {\goth a} = 1$.  However, we prove it in the general case.
Our main application is to the graded case, where we assume that $\goth
a$ and $f$ are homogeneous; then we obtain an isomorphism of graded modules
$Hom_A ({\goth a},A) \cong {\goth b}(deg \ f)$.

Consider the exact sequence
\[
\begin{array}{ccccccccc}
0 & \rightarrow & A & \rightarrow & {\goth a} & \rightarrow & {\goth a} /fA &
\rightarrow & 0 \\
&&1 & \mapsto & f
\end{array}
\]
Since $f$ is $A$-regular, dualizing provides
\[
\begin{array}{ccccccccccc}
0 & \rightarrow & Hom_A ({\goth a} /fA ,A) & \rightarrow & Hom_A ({\goth a} ,A)
& \buildrel \beta \over \rightarrow & Hom_A (A,A) \\
&& || &&&& \phantom{\wr} || \wr \\
&& 0 &&&& A
\end{array}
\]
We first prove that, up to the isomorphism $Hom_A (A,A) \cong A$, we get $Hom_A
({\goth a},A)
\subset {\goth b}$.  Let $\phi \in Hom_A ({\goth a},A)$ and let
$\psi = \beta (\phi)$.  Let $b
:= \psi (1) = \phi(f)$.  Then for any $a \in A$ we have
\[
\psi (a) = \phi(f\cdot a) = a \cdot b.
\]
For any $a \in {\goth a}$ we have
\[
f \cdot \phi(a) = \phi (f\cdot a) = \psi (a) = a \cdot b.
\]
Hence $b \cdot {\goth a} \subset f \cdot A$, i.e.\ $b \in fA :_A {\goth a} =
{\goth b}$.  It follows that $Hom_A ({\goth a},A) \cong im \ \beta \subset
{\goth b}$.

For the reverse inclusion we can define for any $b \in {\goth b}$ a
homomorphism $\phi \in Hom_A ({\goth a},A)$ as the composition of
\[
\begin{array}{ccccccc}
{\goth a} & \rightarrow & fA & \hbox{ and } & fA & \buildrel \sim \over
\rightarrow & A \\
a & \mapsto & ab
\end{array}
\]
Then $\phi (f) = b$.  We conclude that ${\goth b} = im \ \beta \cong Hom_A
({\goth a},A)$.  \qed

\begin{thm} \label{std-and-good}
Suppose that $r+1 \geq 3$.  Then
\newcounter{temp3}
\begin{list} {(\alph{temp3})}{\usecounter{temp3}}
    \setlength{\rightmargin}{\leftmargin}
\item $X$ is  standard determinantal of codimension $r+1$ if and only if there
is a good determinantal subscheme $S \subset \proj{n}$ of codimension $r$ such
that $X \subset S$ is the zero-locus of a regular section  $t \in H^0_* (S,
\widetilde M_S ) = M_S$ for  some $M_S \in {\cal M}_S$.

\item $X$ is good determinantal of codimension $r+1$ if and only if there is a
good determinantal subscheme $S \subset \proj{n}$ of codimension $r$, such
that
$X \subset S$ is the zero-locus of a regular section  $t \in H^0_* (S,
\widetilde M_S ) = M_S$ for some $M_S \in {\cal M}_S$, and the cokernel of this
section is isomorphic to an ideal sheaf in ${\cal O}_X$ of positive height.

\end{list}
\end{thm}

\noindent {\em Proof:}
We first assume that $X$ is standard determinantal and we let
$\Phi$ be a $t \times (t+r)$ homogeneous matrix with $I(\Phi ) = I_X$.  Adding
a general row to $\Phi$ gives a homogeneous $(t+1) \times (t+r)$ matrix
$\Psi$ whose ideal of maximal minors defines a good determinantal scheme $S
\supset X$ of codimension $r$.  We have the commutative diagram
\[
\begin{array}{cccccccccccc}
0 & \rightarrow & ker \ \Psi & \rightarrow & F & \buildrel \Psi \over
\rightarrow & G & \rightarrow & M_S & \rightarrow & 0 \\
&&&& || && \downarrow \\
0 & \rightarrow & ker \ \Phi & \rightarrow & F & \buildrel \Phi \over
\rightarrow & G' & \rightarrow & M_X & \rightarrow & 0 \\
&&&&&&\downarrow \\
&&&&&& 0
\end{array}
\]
where $rk \ F = t+r$, $rk \ G' = t$ and $rk \ G = t+1$.  As in
Theorem~\ref{good-iff-sect}, after possibly twisting we get the exact sequence
\begin{equation}\label{sect-coker}
0 \rightarrow R/I_S (- \hbox{deg } t) \buildrel t \over \rightarrow M_S
\rightarrow M_X
\rightarrow 0.
\end{equation}
Since $S$ is good by construction, Proposition~\ref{good-det} shows that
Lemma~\ref{alg-lemma} applies, setting $A := R/I_S$ and ${\goth a} = M_S$.
This gives
\[
Hom_A (M_S ,A)(-\hbox{deg } t) \cong Ann_A (M_X ) \cong I_X /I_S .
\]
Now, dualizing (\ref{sect-coker}) we get
\[
\begin{array}{cccccccc}
0 & \rightarrow & Hom_A (M_X ,A) & \rightarrow & Hom_A (M_S ,A) &
\buildrel {t^*} \over \rightarrow A (\hbox{deg } t) \\
&& || \\
&& 0
\end{array}
\]
It follows that $X$ is the zero-locus of $t$, proving the direction
$\Rightarrow$ for case (a).  In case (b), we are done by applying
Proposition~\ref{good-det}.

We now consider the direction $\Leftarrow$.  Again let $A = R/I_S$, where $I_S
= I(\Psi)$ for some homogeneous $(t+1) \times (t+r)$ matrix $\Psi$, and apply
the mapping cone construction to the diagram
\[
\begin{array}{cccccccc}
&&&& 0 \\
&&&& \downarrow \\
&& R & \rightarrow & A & \rightarrow & 0 \\
&& \downarrow && \phantom{t} \downarrow t \\
F & \buildrel \Psi \over \rightarrow & G & \rightarrow & M_S \\
&&&& \downarrow \\
&&&& coker \ t \\
&&&& \downarrow \\
&&&& 0
\end{array}
\]
where $rk \ G = t+1$.  This gives the exact sequence
\[
\cdots \rightarrow F \oplus R \buildrel \Phi \over \rightarrow G \rightarrow
coker \ t \rightarrow 0
\]
Since $S$ is good, Proposition~\ref{good-det} gives us that $coker \ t
\cong M_S
/f\cdot A$ for some $A$-regular element $f \in A$ (see the proof of (c)
$\Rightarrow$ (b)).  It follows that $Ann_R (coker \ t)$ has grade $\geq 1 +
grade \ I_S = r+1$, thus $grade \ I(\Phi) = r+1$.  Let $Y$ be the subscheme
defined by $I(\Phi)$.  Then we get as above that $Y$ is the zero-locus of $t$,
and so $X = Y$, and we are done in case (a).  For case (b), again an
application of Proposition~\ref{good-det} completes the argument since $coker \
t \in {\cal M}_X$.  \qed

\bigskip

Note that Theorem~\ref{std-and-good} does not mention global generation, while
Kreuzer's theorem mentioned in the introduction does.
Conjecture~\ref{glob-gen-conj} and Remark~\ref{on-conj} address this.

\begin{conj}\label{glob-gen-conj}
Given $X$ a standard determinantal scheme as in Theorem~\ref{std-and-good}, one
can choose $S$ and $M_S \in {\cal M}_S$ such that $X \subset S$ is the
zero-locus of a regular section  $t \in H^0 (S, \widetilde M_S )$ and such that
$\widetilde M_S$ is globally generated.
\end{conj}

\begin{remark}\label{on-conj}
Consider a free presentation of $M_X$ as in the proof of
Theorem~\ref{std-and-good}:
$$
0 \rightarrow B \rightarrow F \buildrel \Phi \over \rightarrow G
\rightarrow M_X
\rightarrow  0.
$$
Suppose that $\widetilde G$ is globally generated and furthermore that
$\widetilde B^*$ has a regular section $s$.  Then we can write
$$
0 \rightarrow {\cal O} \buildrel s \over \rightarrow \widetilde B^*
 \rightarrow {\cal Q} \rightarrow  0.
$$
A mapping cone gives a
free resolution
$$
0 \rightarrow {\cal O} \oplus \widetilde G^* \rightarrow
\widetilde F^* \rightarrow {\cal Q} \rightarrow  0.
$$
Dualizing this sequence
gives
$$
0 \rightarrow {\cal Q}^* \rightarrow \widetilde F \buildrel \Psi \over
\rightarrow {\cal O}  \oplus \widetilde G \rightarrow {\cal E}xt^1({\cal Q},
{\cal O})
\rightarrow  0.
$$
Since $s$ is a regular section, ${\cal E}xt^1({\cal Q}, {\cal
O})$ is supported on a scheme of codimension one less than the codimension of
$X$. We conclude that $\Psi$ is a Buchsbaum-Rim matrix, and hence $\widetilde
M_S ={\cal E}xt^1({\cal Q}, {\cal O})$ for the scheme $S$ defined by the
maximal
minors of
$\Psi$. As in the proof of Theorem~\ref{std-and-good}, we obtain the exact
sequence
$$
0 \rightarrow R/I_S \rightarrow M_S \rightarrow M_X
\rightarrow 0.
$$
Since ${\cal O} \oplus \widetilde G$ is globally generated, we
see that $\widetilde M_S$ is  globally generated as an ${\cal O}$-module (and
hence as an ${\cal O}_S$-module).

We have just shown that Conjecture~\ref{glob-gen-conj} is true whenever we can
simultaneously guarantee that $\widetilde M_X$ is globally generated and
$\widetilde B^*$ has a regular section. Note in particular that $\widetilde
B^*$
will have a regular section if $\widetilde F^*$ is globally generated.  The
latter holds true, for example, if $X$ is a complete intersection and we choose
$M_X = R/I_X$.
\end{remark}

\begin{remark}\label{CM-type}
Analyzing the proof of Theorem~\ref{std-and-good} and noting that $X$ and $S$
are defined by the maximal minors of a $t \times (t+r)$ matrix and a $(t+1)
\times (t+r)$ matrix, respectively, one observes that there is the following
relation between the Cohen-Macaulay types of $X$ and
$S$, respectively:
\begin{quote} \begin{em}
$X$ has Cohen-Macaulay type ${r + t - 1
  \choose r}$ $\Leftrightarrow$ $S$ has Cohen-Macaulay type ${r + t -1
  \choose r - 1}$.
\end{em}
\end{quote}
This follows from the corresponding Eagon-Northcott resolutions.
\end{remark}


\section{Applications and Examples} \label{corollaries}

In this section we draw some consequences of the results we have shown. We
begin with a characterization of complete intersections. It is well-known
that every complete intersection is arithmetically Gorenstein but the
converse fails in general unless the subscheme has codimension two. For
subschemes of higher codimension we have:

\begin{cor} \label{complete_intersection}
Let $X \subset \proj{n}$ be a subscheme of codimension $r+1 \geq 3$. Then
$X$ is a complete intersection if and only if $X$ is arithmetically
Gorenstein and there is a good determinantal subscheme $S \subset \proj{n}$
of codimension $r$ such
that $X \subset S$ is the zero-locus of a regular section  $t \in H^0_* (S,
\widetilde M_S ) = M_S$ for  some $M_S \in {\cal M}_S$. Furthermore, $S$ and
$M_S$ can be chosen so that $\widetilde M_S$ is globally generated.
\end{cor}

\noindent {\em Proof:} The result follows immediately from
Theorem~\ref{std-and-good}, Remark~\ref{on-conj}, and Remark~\ref{CM-type}.
\qed
\bigskip

Next, we consider subschemes of low codimension. As remarked after
Definition~\ref{Mv}, in the case of codimension two we know
that ${\cal M}_X$ consists of precisely one element (up to isomorphism).

\begin{cor}\label{codim2}
Suppose $X \subset \proj{n}$ ($n \geq 2$) has codimension two.  Then
\newcounter{temp4}
\begin{list} {(\alph{temp4})}{\usecounter{temp4}}
    \setlength{\rightmargin}{\leftmargin}
\item $X$ is standard determinantal if and only if $X$ is arithmetically
Cohen-Macaulay.

\item The following are equivalent:

\newcounter{temp8}
\begin{list} {(\roman{temp8})}{\usecounter{temp8}}
    \setlength{\rightmargin}{\leftmargin}
\item $X$ is good determinantal;

\item $X$ is arithmetically
Cohen-Macaulay and there are an integer $e \in {\Bbb Z}$ and a section $s \in
H^0 (X, \omega_X (e))$ generating $\omega_X (e)$ outside a subscheme of
codimension 3 as an ${\cal O}_X$-module and such that $s$ is a minimal
generator of $H^0_* (\omega_X )$;

\item $X$ is arithmetically Cohen-Macaulay and a generic complete
  intersection .

\end{list}

\end{list}

\end{cor}

\noindent {\em Proof:}
Part (a) is just the Hilbert-Burch theorem.  For (b), the fact that the
codimension of $X$ is 2 implies that $\widetilde M_X \cong \omega_X (e)$ for
some $e \in {\Bbb Z}$.  Then (b) is just a corollary of \linebreak
Proposition~\ref{good-det} and Theorem~\ref{good-iff-glci}.  \qed

\begin{cor}
Suppose that $X \subset \proj{n}$ has codimension 3.  Then $X$ is good
determinantal if and only if there is a good determinantal subscheme $S \subset
\proj{n}$ of codimension 2 such that $X \subset S$ is the zero-locus of a
regular section $t \in H^0 (S, \omega_S (e))$ (for suitable $e \in {\Bbb Z}$)
whose cokernel is supported on a subscheme of codimension $\geq 4$ and
isomorphic to an ideal sheaf of ${\cal O}_X$.
\end{cor}

\noindent {\em Proof:}
This is immediate from Theorem~\ref{std-and-good}.  \qed

\begin{rmk}
\begin{rm}
In general, if $X$ is a good determinantal subscheme of codimension $r+1$ in
$\proj{n}$ then there is a flag of {\em good} determinantal subschemes $X_i$ of
codimension $i$:
\[
X = X_{r+1} \subset X_r \subset \cdots \subset X_2 \subset X_1 \subset
\proj{n}.
\]
In the next corollary we will  show that we can choose the various $X_i$ in
such a way that they have even better properties than guaranteed by the results
of the previous section. \ \ \ \ \ \ \hbox{$\rlap{$\sqcap$}\sqcup$}
\end{rm}
\end{rmk}

\begin{cor}\label{exists-glci-scheme}
If $X \subset \proj{n}$ has codimension $ r + 1 \geq 2$ then the following are
equivalent:
\newcounter{temp9}
\begin{list} {(\alph{temp9})}{\usecounter{temp9}}
    \setlength{\rightmargin}{\leftmargin}
\item $X$ is good determinantal;

\item There is a good determinantal subscheme $S$ of codimension $r$ which
is a local complete
intersection outside a subscheme of codimension $r + 2$, and a section
$t \in H^0_* (S, \widetilde M_S )$ inducing an exact sequence
\[
0 \rightarrow {\cal O}_S(e) \buildrel t \over \rightarrow \widetilde M_S
\rightarrow \widetilde M_X \rightarrow 0
\]
for suitable $M_S \in {\cal M}_S$ and $M_X \in {\cal M}_X$.

\end{list}
\end{cor}

\noindent {\em Proof:}
We first prove (a) $\Rightarrow $ (b).   The existence of a good determinantal
subscheme $S$ and a section $t$ as in the statement follows from
Theorem~\ref{std-and-good} and the exact sequence (\ref{sect-coker}) in
particular.  The only thing remaining to prove is that $S$ can be chosen to
be a
local complete intersection outside a subscheme of codimension $r + 2$ (rather
than codimension $r+1$, as guaranteed by
Proposition~\ref{good-det}).

Assume that the matrix $\Phi$, whose maximal minors define $X$, is a
homogeneous
$t \times (t+r)$ matrix.  The scheme
$S$ is constructed in Theorem~\ref{std-and-good} by adding a ``general row'' to
$\Phi$, producing a $(t+1) \times (t+r)$ matrix, $\Psi$.  One of the points of
the proof of Theorem~\ref{good-iff-glci} is that the locus $Y$ where $S$
fails to
be a local complete intersection is a subscheme of the scheme defined by the
ideal of submaximal minors of $\Psi$.  In particular, $Y$ is a subscheme of
$X$.  The fact that $S$ can be chosen to be a local complete intersection
outside a subscheme of codimension $r+2$ will then follow once we show that,
given a general point $P$ in any component of $X$, there is at least one
submaximal minor of $\Psi$ that does not vanish at $P$.

Since $X$ is good, after a change of basis if necessary we may assume that
there
is a $(t-1) \times (t+r)$ submatrix $\Phi'$ whose ideal of maximal minors
defines a scheme of codimension $r+2$ which is disjoint from $P$.  Hence there
is a maximal minor $A$ of
$\Phi'$ which does not vanish at $P$.  (We make our change of basis, if
necessary, before adding a row to construct
$\Psi$.  Note that we formally include the possibility that $t=1$, i.e.\ that
$X$ is a complete intersection-- see Remark~\ref{one-minor},
Remark~\ref{free-is-BR} and Theorem~\ref{good-iff-sect}.)  Concatenate another
column of $\Phi'$ to
$A$ (by abuse we denote by $A$ both the submatrix and its determinant),
forming a
$(t-1) \times t$ submatrix of $\Phi'$.  Now concatenate the corresponding
elements of the ``general row'' to this matrix, forming a $t \times t$ matrix,
$B$, whose determinant is a submaximal minor of $\Psi$.  Expanding along this
latter row and using the fact that its elements were chosen generally and that
$A$ does not vanish at $P$, we get that the determinant of $B$ does not vanish
at $P$, as desired.  This completes the proof that (a) $\Rightarrow$ (b).

The converse follows exactly as in the proof of
Theorem~\ref{std-and-good} (b).  Note that the condition of being a local
complete intersection away from a subscheme of codimension $r+2$ is  irrelevant
in this direction.
\qed

\begin{rmk}
\begin{rm}
(i) Using the notation of the previous proof we have seen that given a good
determinantal subscheme $X$ we can find  subschemes $Y, S$ such that $Y
\subset X \subset S$ have decreasing codimensions, $X$ is  the zero-locus
of a section of $H^0_* (S, \widetilde M_S )$ and  $X, S$ are local complete
intersections outside $Y$. In this situation we want to call  $X$  a
{\em Cartier divisor on $S$ outside $Y$}. If $Y$ is empty then $X$ is a Cartier
divisor on $S$ in the usual sense.

(ii) Let $X$ be  a good determinantal subscheme of codimension $r+1$ in
$\proj{n}$ and let $X_{r+2} \subset X$ be a subscheme of codimension $r+2$
such that $X$ is a local complete intersection outside $X_{r+2}$. Then
Corollary \ref{exists-glci-scheme} implies that there is a flag of {\em
  good} determinantal subschemes $X_i$ of
codimension $i$:
\[
X = X_{r+1} \subset X_r \subset \cdots \subset X_2 \subset X_1 \subset
X_0 = \proj{n}
\]
such that $X_{i+1}$ is a Cartier divisor on $X_i$ outside $X_{i+2}$ for all
$i = 0,\dots,r$. \ \ \ \ \ \ \hbox{$\rlap{$\sqcap$}\sqcup$}
\end{rm}
\end{rmk}

\begin{cor}\label{exists-lci-curve}
If $X \subset \proj{n}$ is zero-dimensional then the following are equivalent:
\newcounter{temp6}
\begin{list} {(\alph{temp6})}{\usecounter{temp6}}
    \setlength{\rightmargin}{\leftmargin}

\item $X$ is good determinantal;

\item There is a good determinantal curve $S$ which is a local complete
intersection such that $X$ is a Cartier divisor on $S$ associated to  a section
$t \in H^0_* (S, \widetilde M_S )$ inducing an exact sequence
\[
0 \rightarrow {\cal O}_S (e) \buildrel t \over \rightarrow \widetilde M_S
\rightarrow {\cal O}_X (f) \rightarrow 0.
\]

\end{list}
\end{cor}

\noindent {\em Proof:}
Note that under the hypotheses that $X$
is zero-dimensional and good, we get in the commutative diagram
(\ref{usual-diag})  that $coker \
\Phi'$ has finite length, and hence its sheafification is zero.  Hence by the
exact sequence (\ref{short-exact}), we get that the sheafification of $coker \
\Phi$ is just
${\cal O}_X$. Then the result follows from Corollary~\ref{exists-glci-scheme}.
\qed

\bigskip

\begin{cor} \label{points_in_three_space}
Suppose $X \subset \proj{3}$ is zero-dimensional.  Then the following are
equivalent:
\newcounter{temp7}
\begin{list} {(\alph{temp7})}{\usecounter{temp7}}
    \setlength{\rightmargin}{\leftmargin}
\item $X$ is good determinantal;

\item There is an arithmetically Cohen-Macaulay curve $S$, which is a local
  complete
intersection, such that $X$ is a subcanonical Cartier divisor on $S$.

\end{list}
Furthermore, $X$ is defined by a $t \times
(t+r)$ matrix if and only if the Cohen-Macaulay type of $X$ is ${{r+t-1}
\choose
r}$ and that of $S$ is ${{r+t-1}
\choose {r-1}}$.
\end{cor}

\noindent {\em Proof:} Since $S$ has codimension two the exact sequence in
the previous result specializes to the  sequence
\[
0 \rightarrow {\cal O}_S (e) \buildrel t \over \rightarrow \omega_S
\rightarrow {\cal O}_X (f) \rightarrow 0
\]
by Corollary~\ref{codim2}.
Since $S$ is a local complete intersection it implies that $X$ is
subcanonical.  The statement about the Cohen-Macaulay types is just
Remark~\ref{CM-type}.
\qed

\begin{rmk}
\begin{rm}
In view of Remark~\ref{on-conj} and Remark~\ref{CM-type}, Corollaries
\ref{complete_intersection},
\ref{exists-lci-curve} and \ref{points_in_three_space} are generalizations of
Theorem~1.3 of \cite{kreuzer}. \qed
\end{rm}
\end{rmk}

\begin{example}\label{good-ex}
In view of Theorem~\ref{good-iff-glci}, we give examples of curves in
$\proj{3}$ (both of degree 3) to show that a good determinantal scheme need not
be either reduced or a local complete intersection.  For the first,
consider the
curve defined by the matrix
\[
\left [
\begin{array}{ccccc}
x_0 & x_1 & x_2  \\
 0 & x_0 & x_3
\end{array}
\right ]
\]
For the second, consider the curve defined by the matrix
\[
\left [
\begin{array}{cccc}
-x_3 & x_2 & 0 \\
0 & -x_2 & x_1
\end{array}
\right ]
\]
This is the defining matrix for the ``coordinate axes,'' which fail to be a
complete intersection precisely at the ``origin.'' (Recall that in the
definition
of a good determinantal scheme we allowed for the removal of a {\em
generalized}
row.)
\end{example}

\begin{example}
The point of Corollary~\ref{exists-lci-curve} is that given a zero-scheme $X$,
there is so much ``room'' to choose the curve $S$ containing it, that $S$
can be
assumed to be a local complete intersection even at $X$, where one would
normally expect it to have problems.  One naturally can ask if there is so much
room that $S$ can even be taken to be smooth.  The answer is no: for example,
the zeroscheme in $\proj{3}$ defined by the complete intersection $(X_1^2
,X_2^2 ,X_3^2 )$ lies on no smooth curve.  One can ask, though, if there is any
matrix condition analogous to the main result of \cite{chiantini-orecchia}
which guarantees that a ``general'' choice of $S$ will be smooth.
\end{example}

\begin{example}
Any regular section of any twist of the tangent bundle of $\proj{n}$ defines a
good determinantal zero-scheme in $\proj{n}$, by Theorem~\ref{good-iff-sect}.
In fact, it can be shown that if $\cal E$ is any rank $n$ vector bundle on
$\proj{n}$ with $H^i_* (\proj{n} ,{\cal E}) = 0$ for $1 \leq i \leq n-2$, then
any regular section of $\cal E$ defines a good determinantal zero-scheme in
$\proj{n}$.
\end{example}



\begin{thebibliography}{999}

\bibitem{bruns-vetter} W.\ Bruns and U.\ Vetter, ``Determinantal Rings,'' LNM
1327, Springer-Verlag, 1988.

\bibitem{buchs64} D.\ Buchsbaum, {\em A Generalized Koszul complex, I}, Trans.\
Amer.\ Math.\ Soc. {\bf 111} (1964), 183--196.

\bibitem{BE73} D.\ Buchsbaum and D.\ Eisenbud, {\em Remarks on Ideals and
Resolutions}, Symposia Math.\ XI, pp.\ 193-204, Academic Press, London.

\bibitem{BE77} D.\ Buchsbaum and D.\ Eisenbud, {\em What annihilates a
module?}, J.\ Algebra {\bf 47} (1977), 231--243.

\bibitem{BR} D.\ Buchsbaum and D.S.\ Rim, {\em A Generalized Koszul complex,
II. Depth and multiplicity}, Trans.\ Amer.\ Math.\ Soc. {\bf 111} (1964),
197--224.

\bibitem{chiantini-orecchia} L.\ Chiantini and F.\ Orecchia, {\it Plane
Sections
of Arithmetically Normal Curves in} $\proj{3}$, in ``Algebraic Curves and
Projective Geometry, Proceedings (Trento, 1988),'' Lecture Notes in
Mathematics,
vol.\ 1389, Springer--Verlag (1989), 32--42.

\bibitem{EN} J.\ Eagon, D.G.\ Northcott, {\em Ideals defined by matrices, and a
certain complex associated to them}, Proc.\ Royal Soc.\ {\bf 269} (1962),
188--204.

\bibitem{eisenbud} D.\ Eisenbud, ``Commutative Algebra with a view
toward Algebraic Geometry,'' Springer-Verlag, Graduate Texts in Mathematics 150
(1995).

\bibitem{GLP} L.\ Gruson, R.\ Lazarsfeld and C.\ Peskine, {\it On a theorem of
Castelnuovo, and the equations defining space curves}, Invent.\ Math.\ {\bf 72}
(1983), 491--506.

\bibitem{harris} J.\ Harris, ``Algebraic Geometry,'' Springer-Verlag, Graduate
Texts in Mathematics 133 (1992).

\bibitem{H} R.\ Hartshorne, ``Algebraic Geometry,'' Springer-Verlag, Graduate
Texts in Mathematics 52 (1977).

\bibitem{H2} R.\ Hartshorne, {\it Stable reflexive sheaves}, Math.\ Ann.\
  {\bf 254} (1980), 121-176.

\bibitem{kreuzer} M.\ Kreuzer, {\it On 0-dimensional complete intersections},
Math.\ Ann.\ {\bf 292} (1992), 43-58.

\bibitem{MNP} J.\ Migliore, U.\ Nagel and C.\ Peterson, {\em Buchsbaum-Rim
Sheaves}, in preparation.

\bibitem{mig-pet} J.\ Migliore and C.\ Peterson, {\em A Construction of
Codimension Three Arithmetically Gorenstein Subschemes of Projective Space},
to appear in Trans.\ Amer.\ Math.\ Soc.

\end{thebibliography}
\end{document}